%
%

\input lanlmac.tex
\overfullrule=0pt

%



%




\catcode`\@=11

\font\tenmsa=msam10

\font\sevenmsa=msam7

\font\fivemsa=msam5

\font\tenmsb=msbm10

\font\sevenmsb=msbm7

\font\fivemsb=msbm5

\newfam\msafam

\newfam\msbfam

\textfont\msafam=\tenmsa  \scriptfont\msafam=\sevenmsa

  \scriptscriptfont\msafam=\fivemsa

\textfont\msbfam=\tenmsb  \scriptfont\msbfam=\sevenmsb

  \scriptscriptfont\msbfam=\fivemsb

\def\hexnumber@#1{\ifcase#1 0\or1\or2\or3\or4\or5\or6\or7\or8\or9\or

	A\or B\or C\or D\or E\or F\fi }



\font\teneuf=eufm10

\font\seveneuf=eufm7

\font\fiveeuf=eufm5

\newfam\euffam

\textfont\euffam=\teneuf

\scriptfont\euffam=\seveneuf

\scriptscriptfont\euffam=\fiveeuf

\def\frak{\ifmmode\let\next\frak@\else

 \def\next{\Err@{Use \string\frak\space only in math mode}}\fi\next}

\def\goth{\ifmmode\let\next\frak@\else

 \def\next{\Err@{Use \string\goth\space only in math mode}}\fi\next}

\def\frak@#1{{\frak@@{#1}}}

\def\frak@@#1{\fam\euffam#1}


\edef\msa@{\hexnumber@\msafam}

\edef\msb@{\hexnumber@\msbfam}

\mathchardef\boxdot="2\msa@00

\mathchardef\boxplus="2\msa@01

\mathchardef\boxtimes="2\msa@02

\mathchardef\square="0\msa@03

\mathchardef\blacksquare="0\msa@04

\mathchardef\centerdot="2\msa@05

\mathchardef\lozenge="0\msa@06

\mathchardef\blacklozenge="0\msa@07

\mathchardef\circlearrowright="3\msa@08

\mathchardef\circlearrowleft="3\msa@09

\mathchardef\rightleftharpoons="3\msa@0A

\mathchardef\leftrightharpoons="3\msa@0B

\mathchardef\boxminus="2\msa@0C

\mathchardef\Vdash="3\msa@0D

\mathchardef\Vvdash="3\msa@0E

\mathchardef\vDash="3\msa@0F

\mathchardef\twoheadrightarrow="3\msa@10

\mathchardef\twoheadleftarrow="3\msa@11

\mathchardef\leftleftarrows="3\msa@12

\mathchardef\rightrightarrows="3\msa@13

\mathchardef\upuparrows="3\msa@14

\mathchardef\downdownarrows="3\msa@15

\mathchardef\upharpoonright="3\msa@16

\mathchardef\downharpoonright="3\msa@17

\mathchardef\upharpoonleft="3\msa@18

\mathchardef\downharpoonleft="3\msa@19

\mathchardef\rightarrowtail="3\msa@1A

\mathchardef\leftarrowtail="3\msa@1B

\mathchardef\leftrightarrows="3\msa@1C

\mathchardef\rightleftarrows="3\msa@1D

\mathchardef\Lsh="3\msa@1E

\mathchardef\Rsh="3\msa@1F

\mathchardef\rightsquigarrow="3\msa@20

\mathchardef\leftrightsquigarrow="3\msa@21

\mathchardef\looparrowleft="3\msa@22

\mathchardef\looparrowright="3\msa@23

\mathchardef\circeq="3\msa@24

\mathchardef\succsim="3\msa@25

\mathchardef\gtrsim="3\msa@26

\mathchardef\gtrapprox="3\msa@27

\mathchardef\multimap="3\msa@28

\mathchardef\therefore="3\msa@29

\mathchardef\because="3\msa@2A

\mathchardef\doteqdot="3\msa@2B

\mathchardef\triangleq="3\msa@2C

\mathchardef\precsim="3\msa@2D

\mathchardef\lesssim="3\msa@2E

\mathchardef\lessapprox="3\msa@2F

\mathchardef\eqslantless="3\msa@30

\mathchardef\eqslantgtr="3\msa@31

\mathchardef\curlyeqprec="3\msa@32

\mathchardef\curlyeqsucc="3\msa@33

\mathchardef\preccurlyeq="3\msa@34

\mathchardef\leqq="3\msa@35

\mathchardef\leqslant="3\msa@36

\mathchardef\lessgtr="3\msa@37

\mathchardef\backprime="0\msa@38

\mathchardef\risingdotseq="3\msa@3A

\mathchardef\fallingdotseq="3\msa@3B

\mathchardef\succcurlyeq="3\msa@3C

\mathchardef\geqq="3\msa@3D

\mathchardef\geqslant="3\msa@3E

\mathchardef\gtrless="3\msa@3F

\mathchardef\sqsubset="3\msa@40

\mathchardef\sqsupset="3\msa@41

\mathchardef\vartriangleright="3\msa@42

\mathchardef\vartriangleleft="3\msa@43

\mathchardef\trianglerighteq="3\msa@44

\mathchardef\trianglelefteq="3\msa@45

\mathchardef\bigstar="0\msa@46

\mathchardef\between="3\msa@47

\mathchardef\blacktriangledown="0\msa@48

\mathchardef\blacktriangleright="3\msa@49

\mathchardef\blacktriangleleft="3\msa@4A

\mathchardef\vartriangle="0\msa@4D

\mathchardef\blacktriangle="0\msa@4E

\mathchardef\triangledown="0\msa@4F

\mathchardef\eqcirc="3\msa@50

\mathchardef\lesseqgtr="3\msa@51

\mathchardef\gtreqless="3\msa@52

\mathchardef\lesseqqgtr="3\msa@53

\mathchardef\gtreqqless="3\msa@54

\mathchardef\Rrightarrow="3\msa@56

\mathchardef\Lleftarrow="3\msa@57

\mathchardef\veebar="2\msa@59

\mathchardef\barwedge="2\msa@5A

\mathchardef\doublebarwedge="2\msa@5B

\mathchardef\angle="0\msa@5C

\mathchardef\measuredangle="0\msa@5D

\mathchardef\sphericalangle="0\msa@5E

\mathchardef\varpropto="3\msa@5F

\mathchardef\smallsmile="3\msa@60

\mathchardef\smallfrown="3\msa@61

\mathchardef\Subset="3\msa@62

\mathchardef\Supset="3\msa@63

\mathchardef\Cup="2\msa@64

\mathchardef\Cap="2\msa@65

\mathchardef\curlywedge="2\msa@66

\mathchardef\curlyvee="2\msa@67

\mathchardef\leftthreetimes="2\msa@68

\mathchardef\rightthreetimes="2\msa@69

\mathchardef\subseteqq="3\msa@6A

\mathchardef\supseteqq="3\msa@6B

\mathchardef\bumpeq="3\msa@6C

\mathchardef\Bumpeq="3\msa@6D

\mathchardef\lll="3\msa@6E

\mathchardef\ggg="3\msa@6F

\mathchardef\circledS="0\msa@73

\mathchardef\pitchfork="3\msa@74

\mathchardef\dotplus="2\msa@75

\mathchardef\backsim="3\msa@76

\mathchardef\backsimeq="3\msa@77

\mathchardef\complement="0\msa@7B

\mathchardef\intercal="2\msa@7C

\mathchardef\circledcirc="2\msa@7D

\mathchardef\circledast="2\msa@7E

\mathchardef\circleddash="2\msa@7F

\def\ulcorner{\delimiter"4\msa@70\msa@70 }

\def\urcorner{\delimiter"5\msa@71\msa@71 }

\def\llcorner{\delimiter"4\msa@78\msa@78 }

\def\lrcorner{\delimiter"5\msa@79\msa@79 }

\def\yen{\mathhexbox\msa@55 }

\def\checkmark{\mathhexbox\msa@58 }

\def\circledR{\mathhexbox\msa@72 }

\def\maltese{\mathhexbox\msa@7A }

\mathchardef\lvertneqq="3\msb@00

\mathchardef\gvertneqq="3\msb@01

\mathchardef\nleq="3\msb@02

\mathchardef\ngeq="3\msb@03

\mathchardef\nless="3\msb@04

\mathchardef\ngtr="3\msb@05

\mathchardef\nprec="3\msb@06

\mathchardef\nsucc="3\msb@07

\mathchardef\lneqq="3\msb@08

\mathchardef\gneqq="3\msb@09

\mathchardef\nleqslant="3\msb@0A

\mathchardef\ngeqslant="3\msb@0B

\mathchardef\lneq="3\msb@0C

\mathchardef\gneq="3\msb@0D

\mathchardef\npreceq="3\msb@0E

\mathchardef\nsucceq="3\msb@0F

\mathchardef\precnsim="3\msb@10

\mathchardef\succnsim="3\msb@11

\mathchardef\lnsim="3\msb@12

\mathchardef\gnsim="3\msb@13

\mathchardef\nleqq="3\msb@14

\mathchardef\ngeqq="3\msb@15

\mathchardef\precneqq="3\msb@16

\mathchardef\succneqq="3\msb@17

\mathchardef\precnapprox="3\msb@18

\mathchardef\succnapprox="3\msb@19

\mathchardef\lnapprox="3\msb@1A

\mathchardef\gnapprox="3\msb@1B

\mathchardef\nsim="3\msb@1C


\mathchardef\ncong="3\msb@1D

\mathchardef\varsubsetneq="3\msb@20

\mathchardef\varsupsetneq="3\msb@21

\mathchardef\nsubseteqq="3\msb@22

\mathchardef\nsupseteqq="3\msb@23

\mathchardef\subsetneqq="3\msb@24

\mathchardef\supsetneqq="3\msb@25

\mathchardef\varsubsetneqq="3\msb@26

\mathchardef\varsupsetneqq="3\msb@27

\mathchardef\subsetneq="3\msb@28

\mathchardef\supsetneq="3\msb@29

\mathchardef\nsubseteq="3\msb@2A

\mathchardef\nsupseteq="3\msb@2B

\mathchardef\nparallel="3\msb@2C

\mathchardef\nmid="3\msb@2D

\mathchardef\nshortmid="3\msb@2E

\mathchardef\nshortparallel="3\msb@2F

\mathchardef\nvdash="3\msb@30

\mathchardef\nVdash="3\msb@31

\mathchardef\nvDash="3\msb@32

\mathchardef\nVDash="3\msb@33

\mathchardef\ntrianglerighteq="3\msb@34

\mathchardef\ntrianglelefteq="3\msb@35

\mathchardef\ntriangleleft="3\msb@36

\mathchardef\ntriangleright="3\msb@37

\mathchardef\nleftarrow="3\msb@38

\mathchardef\nrightarrow="3\msb@39

\mathchardef\nLeftarrow="3\msb@3A

\mathchardef\nRightarrow="3\msb@3B

\mathchardef\nLeftrightarrow="3\msb@3C

\mathchardef\nleftrightarrow="3\msb@3D

\mathchardef\divideontimes="2\msb@3E

\mathchardef\varnothing="0\msb@3F

\mathchardef\nexists="0\msb@40

\mathchardef\mho="0\msb@66

\mathchardef\eth="0\msb@67

\mathchardef\eqsim="3\msb@68

\mathchardef\beth="0\msb@69

\mathchardef\gimel="0\msb@6A

\mathchardef\daleth="0\msb@6B

\mathchardef\lessdot="3\msb@6C

\mathchardef\gtrdot="3\msb@6D

\mathchardef\ltimes="2\msb@6E

\mathchardef\rtimes="2\msb@6F

\mathchardef\shortmid="3\msb@70

\mathchardef\shortparallel="3\msb@71

\mathchardef\smallsetminus="2\msb@72

\mathchardef\thicksim="3\msb@73

\mathchardef\thickapprox="3\msb@74

\mathchardef\approxeq="3\msb@75

\mathchardef\succapprox="3\msb@76

\mathchardef\precapprox="3\msb@77

\mathchardef\curvearrowleft="3\msb@78

\mathchardef\curvearrowright="3\msb@79

\mathchardef\digamma="0\msb@7A

\mathchardef\varkappa="0\msb@7B

\mathchardef\hslash="0\msb@7D

\mathchardef\hbar="0\msb@7E

\mathchardef\backepsilon="3\msb@7F

\def\Bbb{\ifmmode\let\next\Bbb@\else

 \def\next{\errmessage{Use \string\Bbb\space only in math mode}}\fi\next}

\def\Bbb@#1{{\Bbb@@{#1}}}

\def\Bbb@@#1{\fam\msbfam#1}

\catcode`\@=12

\def\Pr{[\hyperref {}{reference}{1}{1}]}
\def\RS{[\hyperref {}{reference}{2}{2}]}
\def\JBZ{[\hyperref {}{reference}{3}{3}]}
\def\macmah{(\relax \hyperref {}{equation}{1.1}{1.1})}
\def\dG{[\hyperref {}{reference}{4}{4}]}
\def\generic{1}
\def\polygon{2}
\def\bijec{3}
\def\degier{4}
\def\fixeda{5}
\def\fixedb{6}
\def\K{[\hyperref {}{reference}{5}{5}]}
\def\fixedc{7}
\def\split{8}
\def\empty{9}
\def\Wie{[\hyperref {}{reference}{6}{6}]}
\def\fulempty{10}
\writedefs
%
\input epsf
\def\fig#1#2#3{%
\xdef#1{\the\figno}%
\writedef{#1\leftbracket \the\figno}%
\nobreak%
\par\begingroup\parindent=0pt\leftskip=1cm\rightskip=1cm\parindent=0pt%
\baselineskip=11pt%
\midinsert%
\centerline{#3}%
\vskip 12pt%
{\bf Fig.\ \the\figno:} #2\par%
\endinsert\endgroup\par%
\goodbreak%
\global\advance\figno by1%
}
\def\omit#1{}

\def\rem#1{{\sl [#1]}}
\def\pre#1{ ({\tt #1})}

\lref\Pr{J. Propp, {\sl The many faces of alternating-sign matrices}, 
preprint\pre{math.CO/0208125}.}
\lref\RS{A.V. Razumov and Yu.G. Stroganov, 
{\sl Combinatorial nature
of ground state vector of O(1) loop model},
preprint\pre{math.CO/0104216}.}
\lref\JBZ{J.-B.~Zuber, {\sl On the Counting of Fully Packed Loop Configurations. Some new conjectures}, preprint\pre{math-ph/0309057}.}
\lref\dG{J. de Gier, {\sl The art of number guessing: where
combinatorics meets physics}, preprint\pre{math.CO/0211285}.}
\lref\K{C. Krattenthaler, private communication.}
\lref\Wie{B. Wieland, {\sl  A large dihedral symmetry of the set of
alternating-sign matrices}, 
 {\it Electron. J. Combin.} {\bf 7} (2000) R37, 
preprint\pre{math.CO/0006234}.}
%
%
\def\ppmove{$\epsfbox{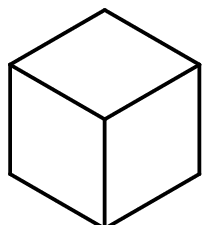}\leftrightarrow\epsfbox{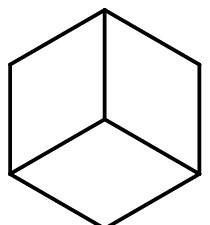}$}
\def\dimermove{$\epsfbox{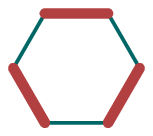}\leftrightarrow\epsfbox{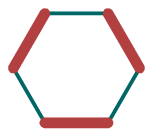}$}
\def\fplmovev{$\epsfbox{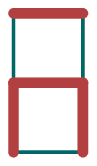}\leftrightarrow\epsfbox{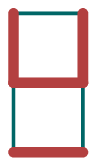}$}
\def\fplmoveh{$\epsfbox{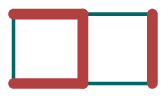}\leftrightarrow\epsfbox{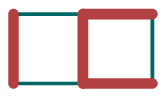}$}

\Title{}
{\vbox{
\centerline{A Bijection between classes of Fully Packed Loops}
\medskip
\centerline{and Plane Partitions}
}}
\bigskip
\centerline{P.~Di~Francesco,\footnote{${}^\#$}{
Service de Physique Th\'eorique de Saclay, 
CEA/DSM/SPhT, URA 2306 du CNRS, 
C.E.A.-Saclay, F-91191 Gif sur Yvette Cedex, France}}
\medskip
\centerline{P. Zinn-Justin \footnote{${}^\star$}
{Laboratoire de Physique Th\'eorique et Mod\`eles Statistiques, UMR 8626 du CNRS,
Universit\'e Paris-Sud, B\^atiment 100,  F-91405 Orsay Cedex, France
}}
\medskip
\centerline{and J.-B. Zuber ${}^\#$}\medskip
\bigskip
\noindent
It has recently been observed empirically that 
the number of FPL configurations with 3 sets of $a$, $b$ and $c$ nested arches 
equals  the number of plane partitions in a box
of size $a\times b \times c$. In this note,  
this result is proved by constructing explicitly the bijection 
between these FPL and plane partitions.


\omit{reste \`a faire: unifier les denominations 
links $\leftrightarrow$ edges $\surd$; 
unoccupied (edges) $\leftrightarrow$ empty  $\surd$;
Choix (voire redessin de certaines?) des figures $\surd$;
References a completer  $\surd$; Ackowledgements $\surd$;
Choix du journal : Quid de Electr. J. Combin? 
C'est la que Wieland a publie le sien. 
Mention du site de PZJ?$\surd$ }
\Date{11/2003}\def\rem#1{}

\newsec{Introduction}

Configurations of Fully Packed Loops, or FPL in short, are 
sets of disconnected paths visiting once each vertex of a 
square $n\times n$ grid, and exiting through every other of the 
$4n$ external edges. While these FPL constitute an interesting and 
much studied model of statistical mechanics, they have also attracted
recently  the attention of combinatorialists by their connections
with alternating sign matrices,  tiling problems, plane partitions and
related topics (see for example \Pr\ for an overall review).
Moreover, there is a yet mysterious relation 
with a linear problem: the numbers of FPL configurations of
different ``link patterns'' give the components of the 
Perron-Frobenius eigenvector of an operator (Hamiltonian) constructed
in terms of the generators of the Temperley-Lieb algebra \RS. 

In a recent work \JBZ, one of us has observed empirically that 
the number of FPL with 3 sets of $a$, $b$ and $c$ nested arches, 
($a+b+c=n$), 
is nothing else than the number of plane partitions in a box
of size $a\times b \times c$  (MacMahon formula)
\eqnn\macmah
$$\eqalignno{PP(a,b,c)
&=
 {{ n-1\choose a}{ n-2\choose a}\cdots{ n-b\choose a}\over
{\rm same \ for \ } n=a+b}\cr
&=
\prod_{i=1}^a\prod_{j=1}^b\prod_{k=1}^c {i+j+k-1\over i+j+k-2}\ \cr
&=
{H(a+b+c) H(a)H(b)H(c)\over
H(a+b)H(b+c-1)H(c+a) }\ .\cr
}$$
In the third of these equivalent expressions, 
$H$ is the hyperfactorial function $H(p)=(p-1)!(p-2)!\cdots 1!$. 
It is the object of this note to
prove this result by constructing explicitly the bijection 
between these FPL and plane partitions.
In fact a similar bijection between FPL with different boundary 
conditions and a tiling problem had been
constructed by de Gier \dG, and our construction is closely
related to his method. 

We now state precisely the result. 
Let us consider a FPL with three sets of nested arches;
let $a$, $b$, $c$ be the numbers of nested arches; and $A$,
$B$ and $C$ be the centers of the three  bundles of arches, 
namely the central unoccupied external edges.
We call such a FPL a ``FPL of type $(A,B,C)$''.
The number of occupied external edges between  $A$ and $B$ is $a+b$. Thus
the data of $A$, $B$, $C$ determine $a$, $b$, $c$.

\font\sc=cmcsc10
\noindent{\sc Theorem.} {\sl There is a bijection between FPL configurations 
of type $(A,B,C)$ (with a link pattern made of three nested 
sets of $a$, $b$ and $c$ arches) 
and the plane partitions in a box of size $a\times b \times c$.}

\newsec{The bijection}
In this section, we explicitly construct
the bijection between FPL configurations of type $(A,B,C)$
and plane partitions in a box of sides $a$, $b$, $c$.

Given the center $A$ of a bundle of  arches, 
on the side of the $n\times n$ square, 
we construct the {\it cone}\/ which is the space between
the two diagonals at 45 degrees starting from the innermost endpoint of the
unoccupied edge $A$. 
One then proves easily that there exists a {\it unique}\/ 
point among $A$, $B$, $C$  such that its cone contains the two others;
this point is never in a corner. 
Up to a permutation of the letters 
$A$, $B$, $C$ (and of the associated $a$, $b$, $c$),
we may always assume that this point is $C$ and then, 
three cases may occur
(see Fig.~\generic\ for the generic case, and \polygon\ 
for a sample of typical cases that we shall follow throughout this paper):

{
\fig\generic{The 
three generic cases, with $A$ and $B$ in the cone of $C$: 
($i$) $A,B$ on the opposite side of $C$; 
($ii$) $A$ on the opposite side of $C$ and $B$ on an adjacent one; 
($iii$) $A$ and $B$ on distinct adjacent sides to $C$. 
In each case, we have represented the polygon $P'$ in thick black
lines.  The distances shown on the figure equal  
the number of elementary segments, possibly 
rounded to one of the nearest integers.
}{\epsfbox{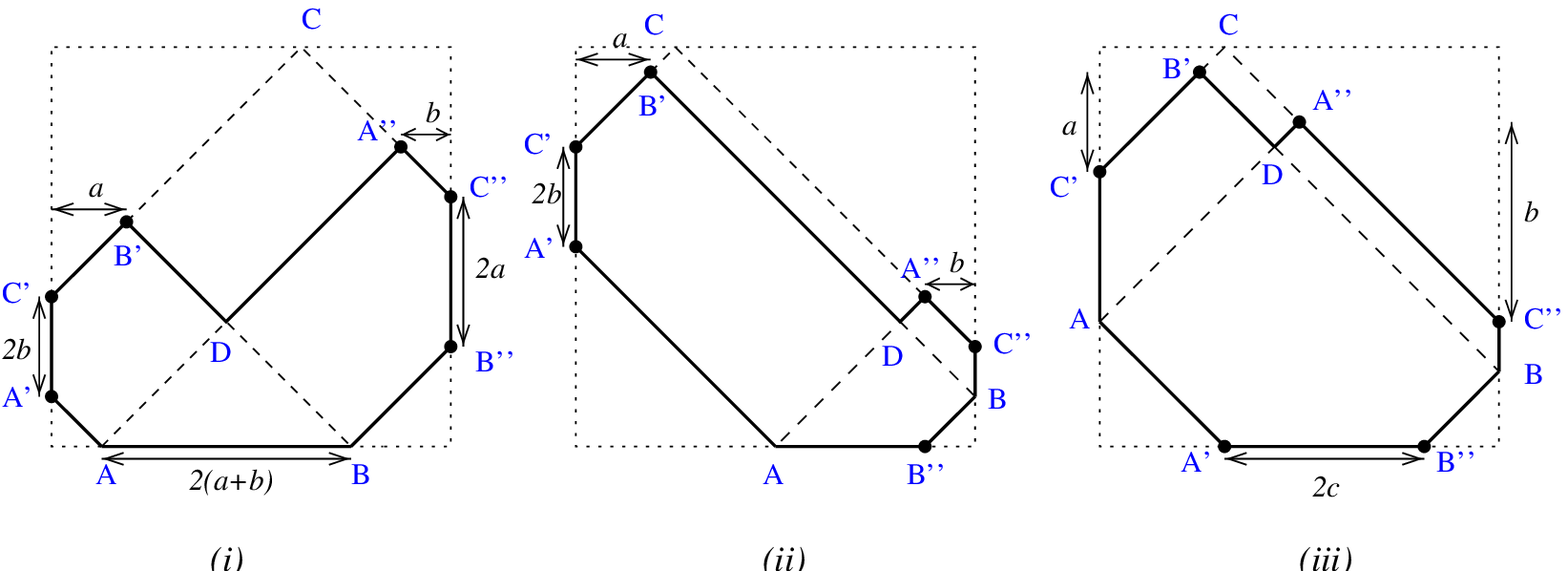}}}

{\item{$(i)$} $A$ and $B$ lie on the same side of the square, which is 
 necessarily  opposite to $C$, since otherwise $a+c$ or  $b+c \ge  n$, 
which is absurd;
\item{$(ii)$} $A$ lies on the side opposite to $C$, and $B$
on an adjacent side: 
any such configuration or its mirror image is as depicted 
on Fig.~\generic~$(ii)$;
\item{$(iii)$} $A$ and $B$ are on sides adjacent to that of $C$, 
which are necessarily opposite.

\noindent Remark: the limiting cases where $A$ and/or $B$ 
is in a corner offer no difficulty and are treated 
as in the general case.

In all that follows, we shall treat in parallel these 3 cases; even
though the reasoning is essentially the same, there are some technical
differences between the 3 situations.

We now describe the 4 steps required to produce the bijection, 
postponing to the next section the actual proof of the theorem.

\fig\polygon{Polygons, dominos and hexagons.
The values of $(a,b,c)$ are respectively
$(2,3,10)$, $(3,2,7)$ and $ (3,4,3)$.  
}{
$\matrix{
\epsfbox{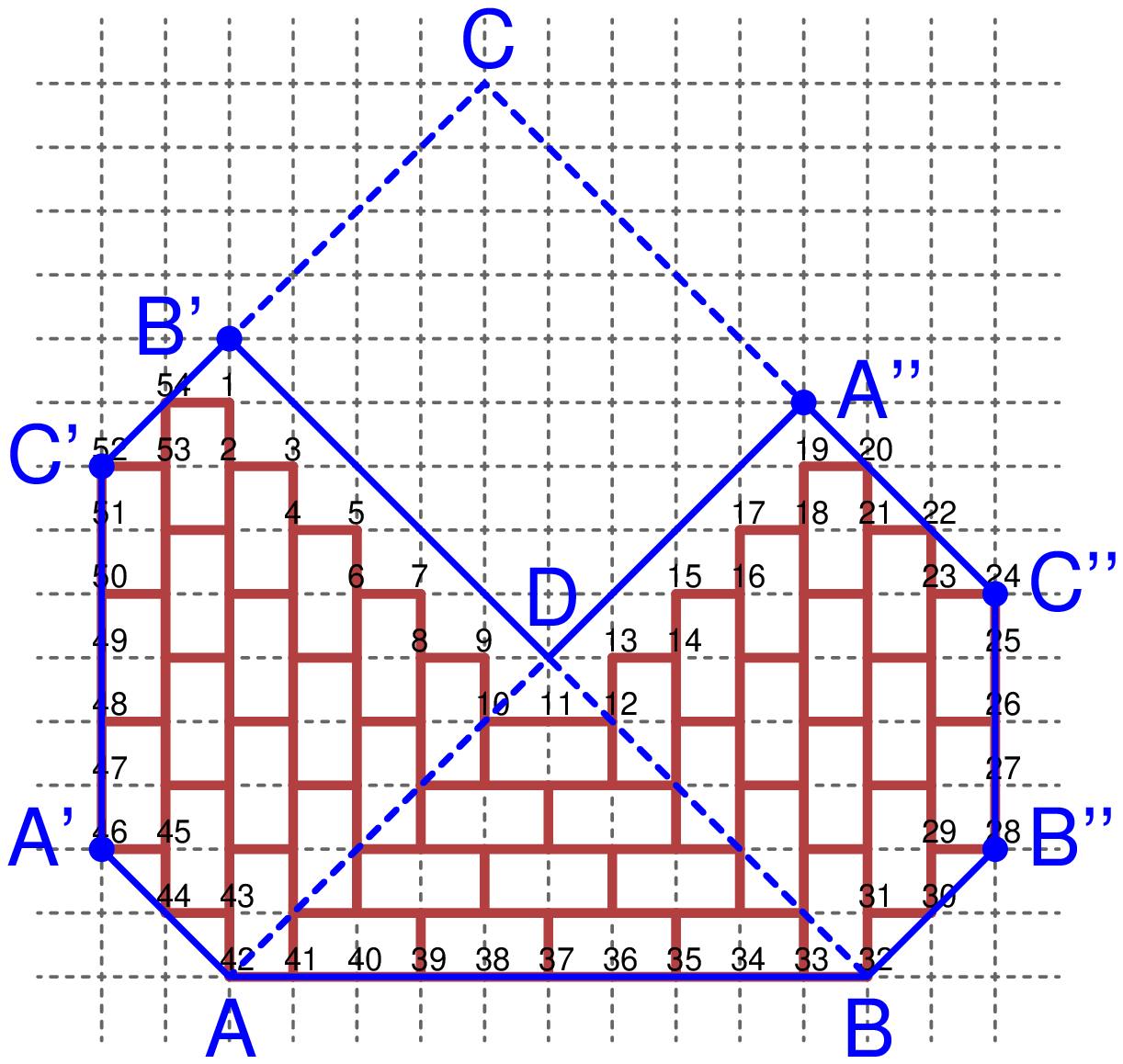}&\epsfbox{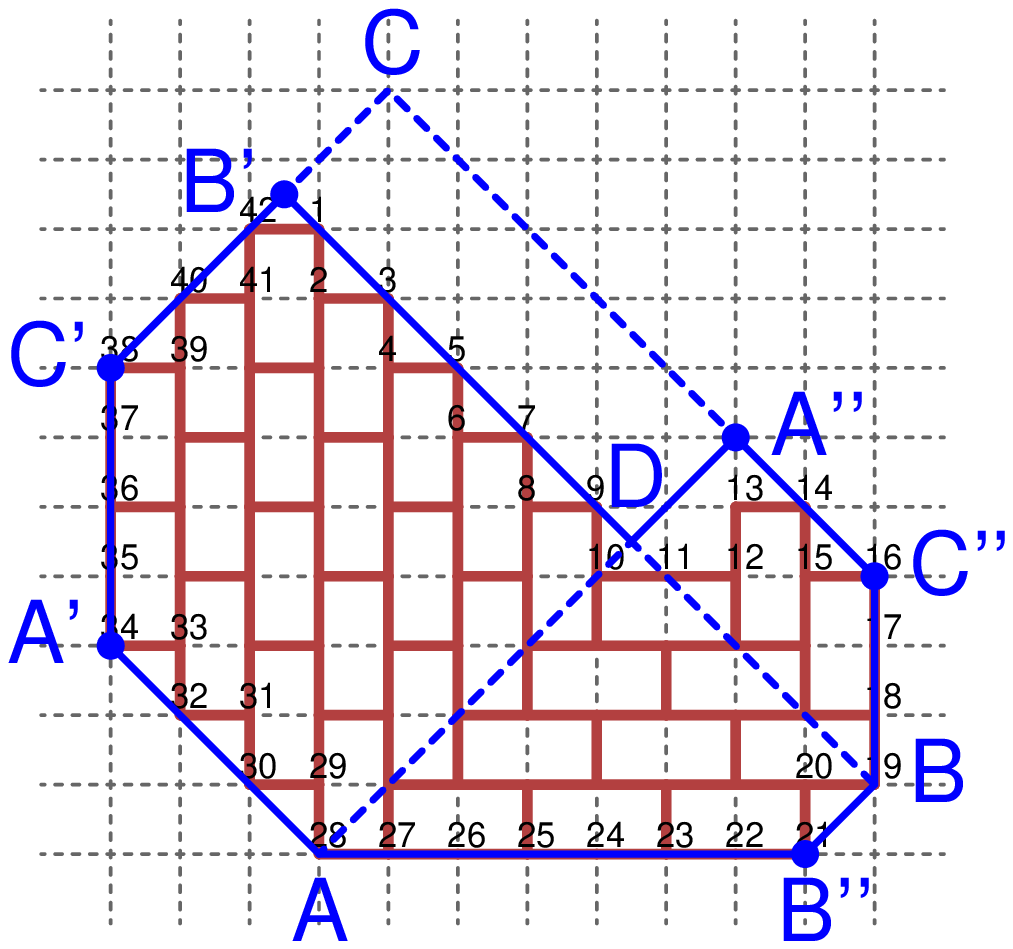}&\epsfbox{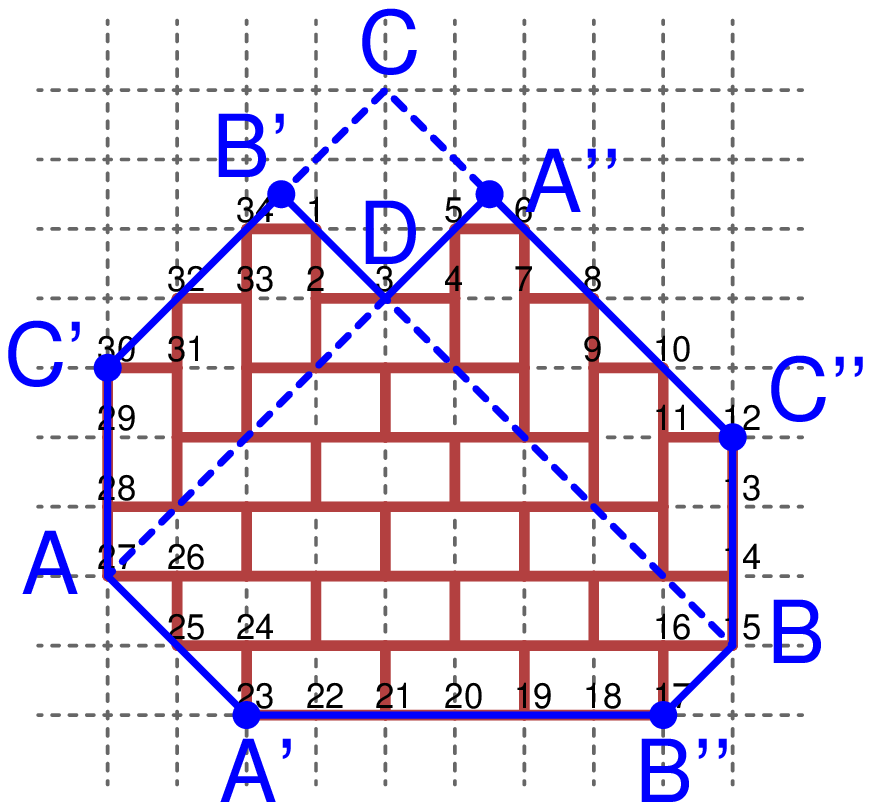}\cr
\epsfbox{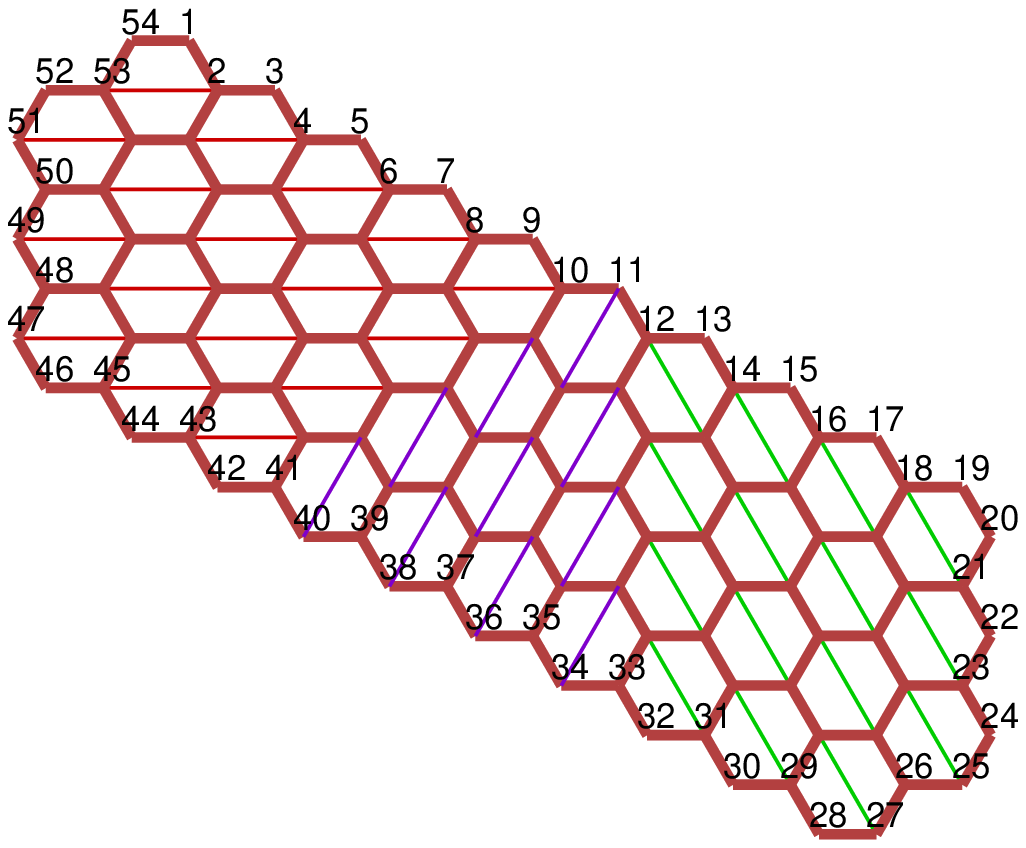}&\epsfbox{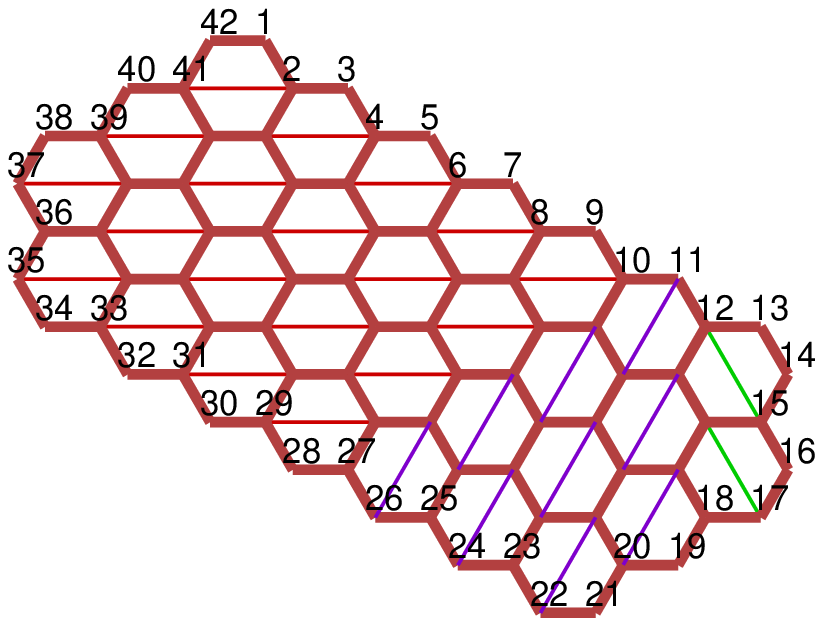}&\epsfbox{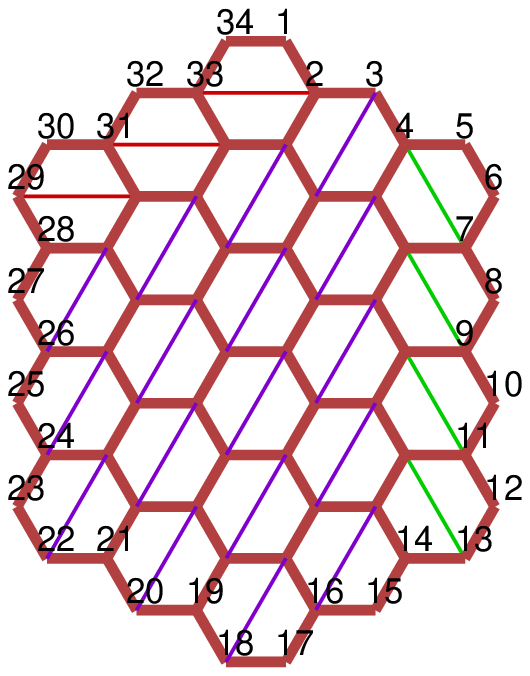}\cr
(i)&(ii)&(iii)\cr
}$ }

\item{\sl Step 1} From the points $C$, draw the diagonals (lines of
slope $\pm 1$; they are the boundaries of the cone mentioned above), 
which cross the sides of the square at $C'$, $C''$.
{}From $A$, resp.\ $B$,  draw also the diagonals, and call $A'$ and $A''$, 
resp.\ $B'$ and $B''$, their intersection with a side of the square 
or with the diagonals coming from $C$, whichever comes first. 
(Some of these points may coincide). Finally let $D$ be 
the intersection of $AA''$ and $BB'$.
Let $P$ be the polygon
$(AA')C'CC''(B''B)$, where the brackets $(AA')$ mean a pair up to
transposition, depending on the case (see Fig.~\generic\ and 
\polygon); and $P'$ be   
the polygon $(AA')C'B'DA''C''(B''B)$.

\item{\sl Step 2} Pave the inside of the polygon $P'$ with
dominos, in the way indicated on Fig.~\polygon\ 
(note that in all cases the polygon
is not quite filled with the dominos due to diagonal lines and in the vicinity of $D$).

\item{\sl Step 3} One can deform the dominos into a subset of the hexagonal
lattice. One can check, see Fig.~\polygon, that in the three
 cases this subset has the desired shape
of a hexagon with sides of lengths $a$, $b$, $c$. 
The middle edges of the dominos become diameters of the hexagons.

\item{\sl Step 4} Finally, to each FPL of type $(A,B,C)$, associate a
dimer configuration
on the subset of the hexagonal lattice by keeping only
the FPL edges of the borders of the dominos, while discarding
all other edges including the middle edges of the dominos;  
and finally deform the dominos into the hexagons as explained in step 3. 
The dimer configuration can equivalently be represented as
a plane partition in a box $a\times b \times c$, see Fig.~\bijec. 
 
\fig\bijec{The bijection: examples.}{
$\matrix{
\epsfbox{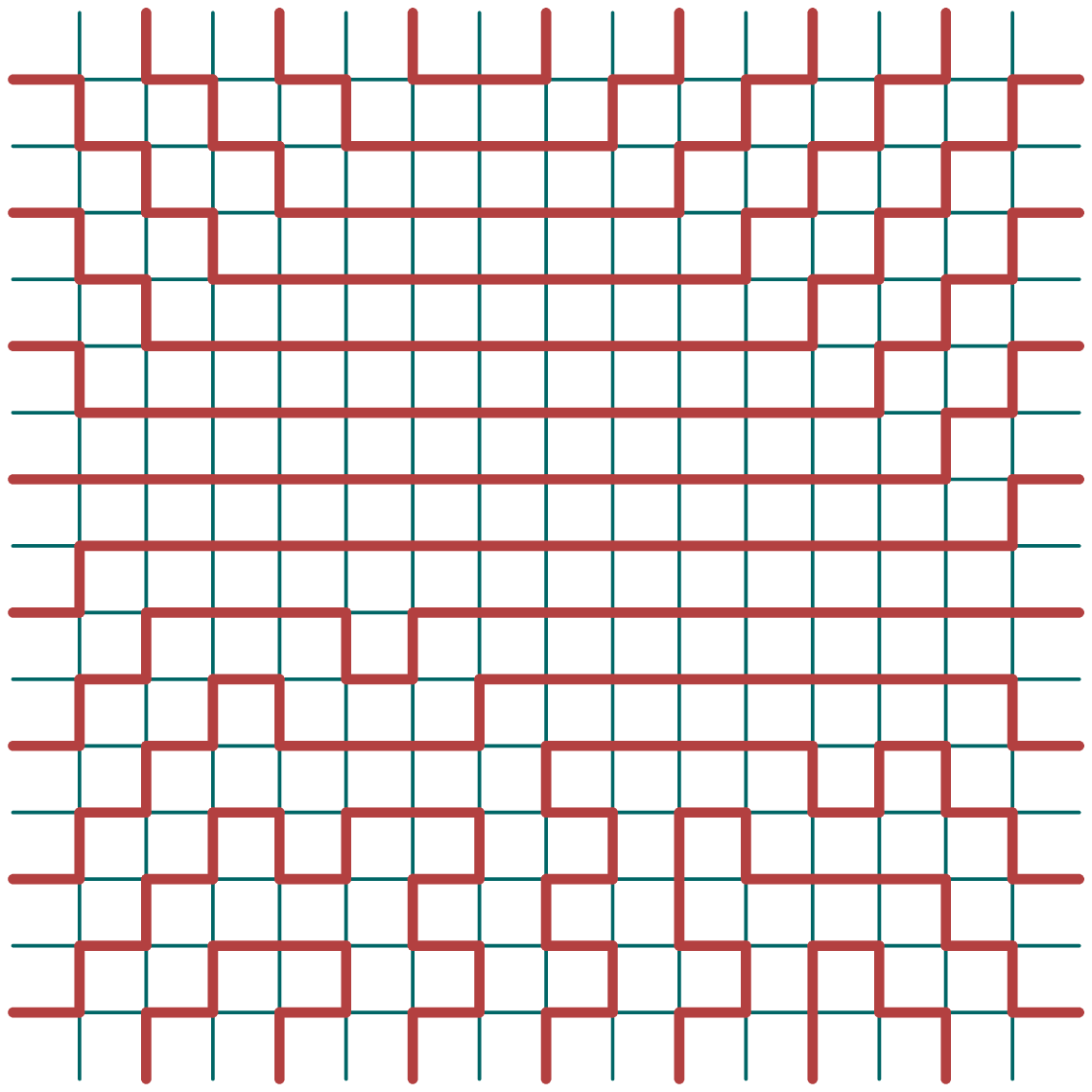}&\epsfbox{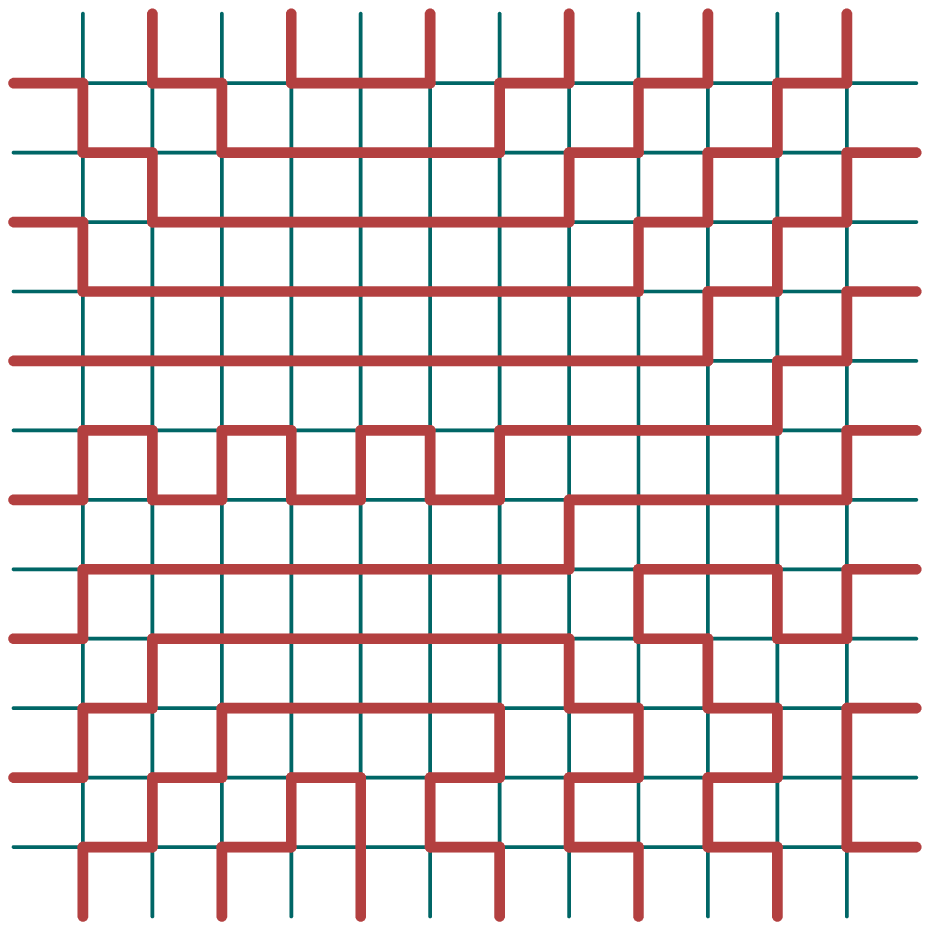}&\epsfbox{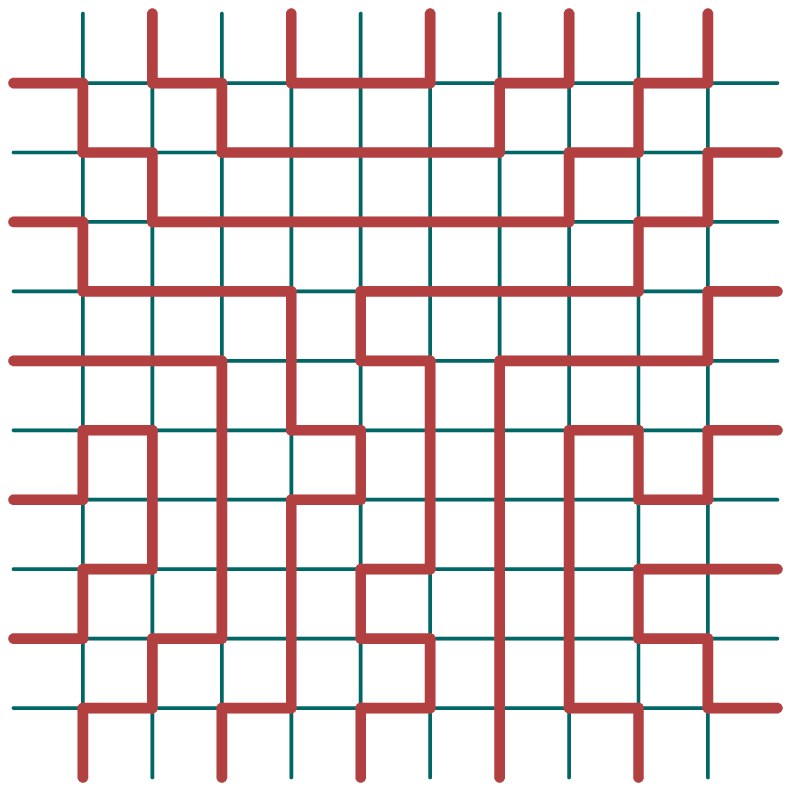}\cr
\epsfbox{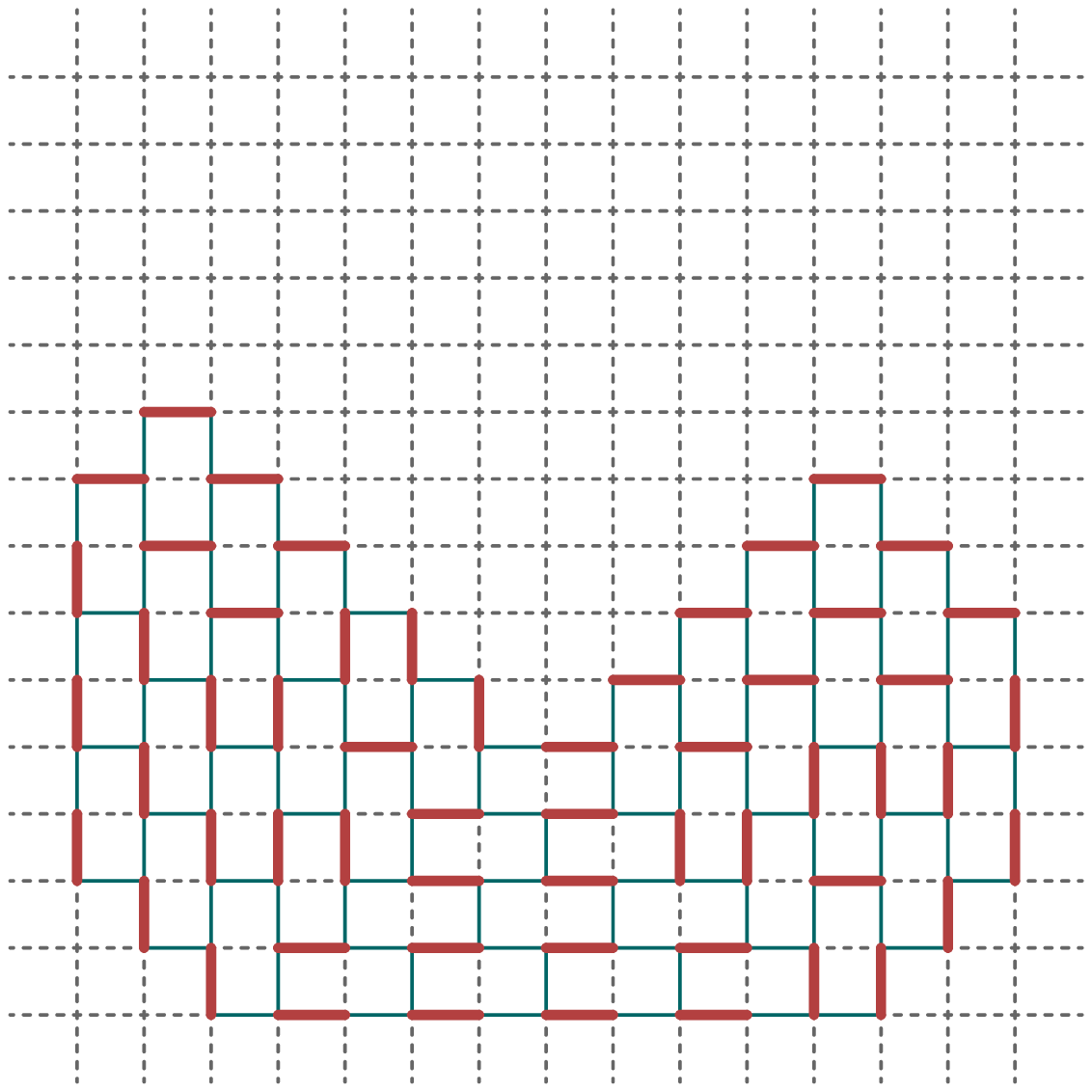}&\epsfbox{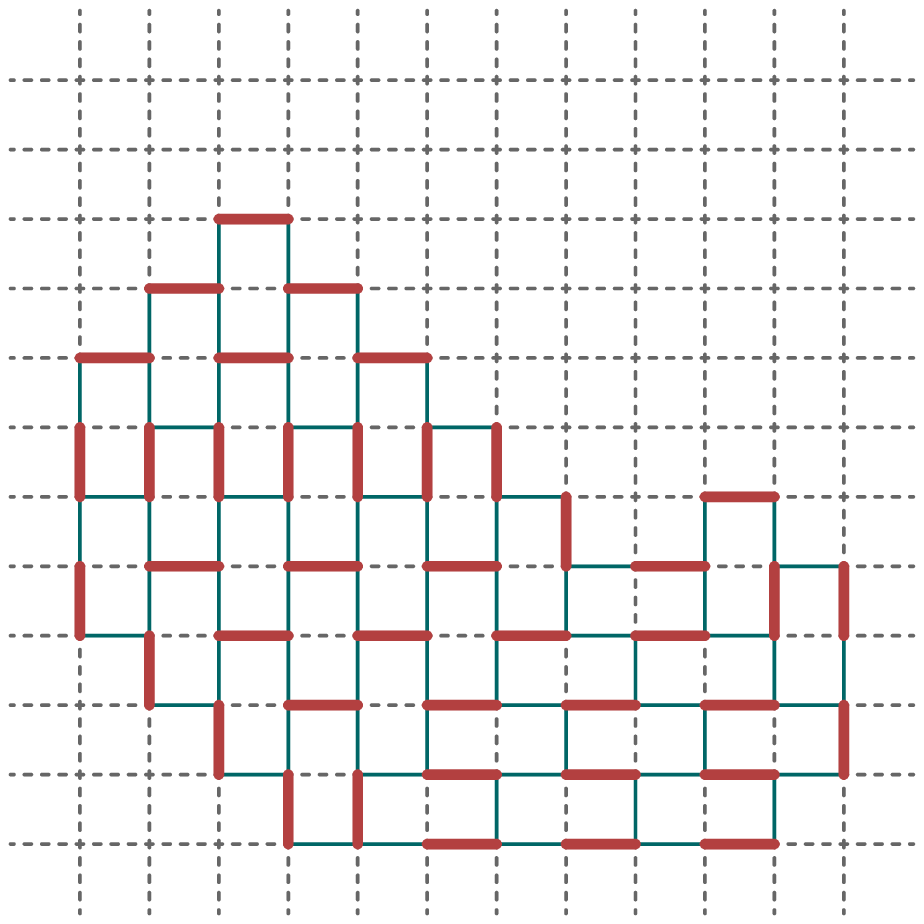}&\epsfbox{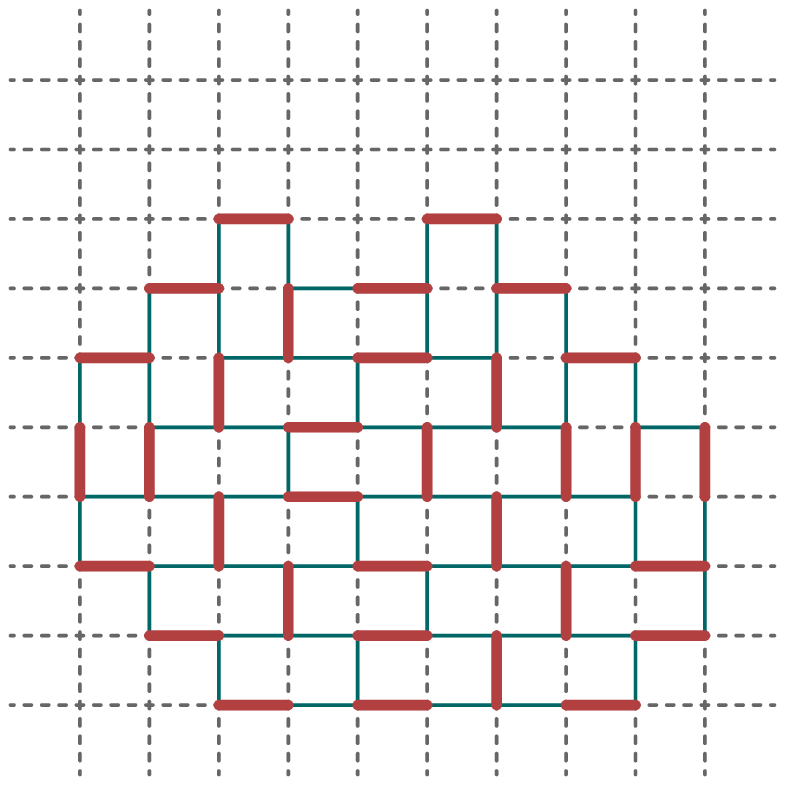}\cr
\epsfbox{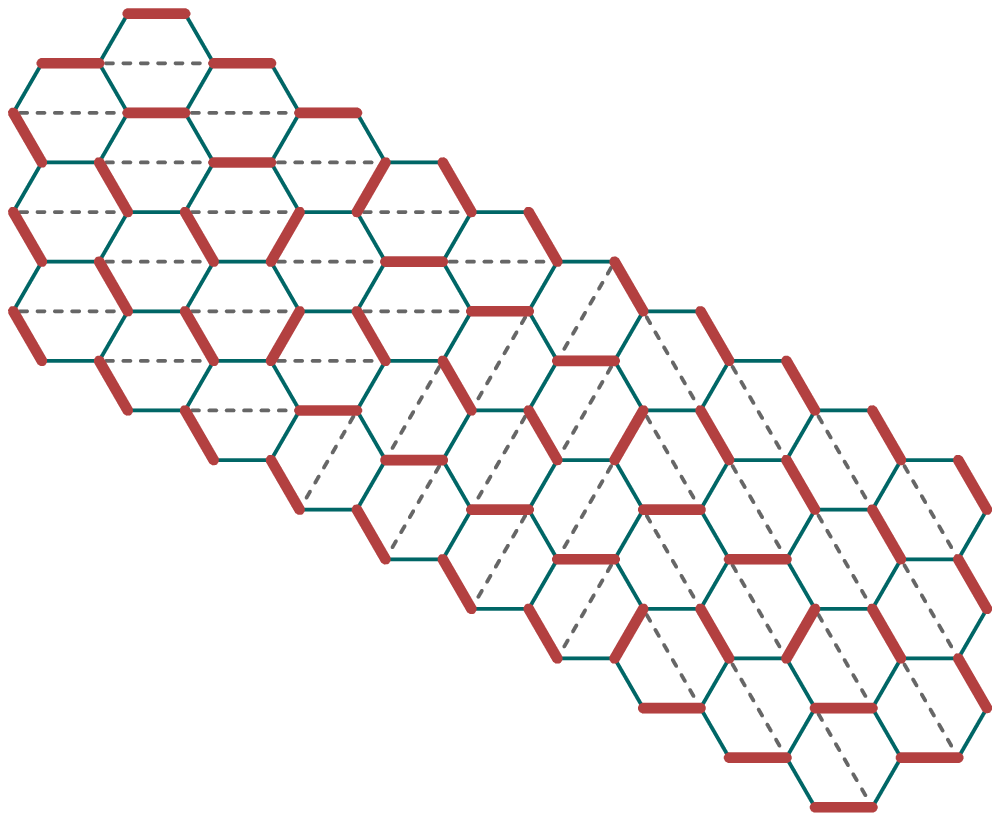}&\epsfbox{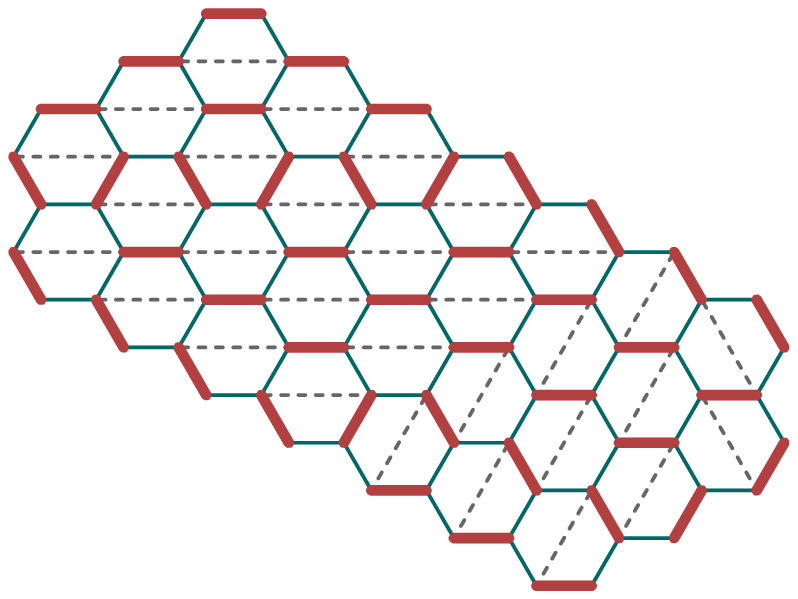}&\epsfbox{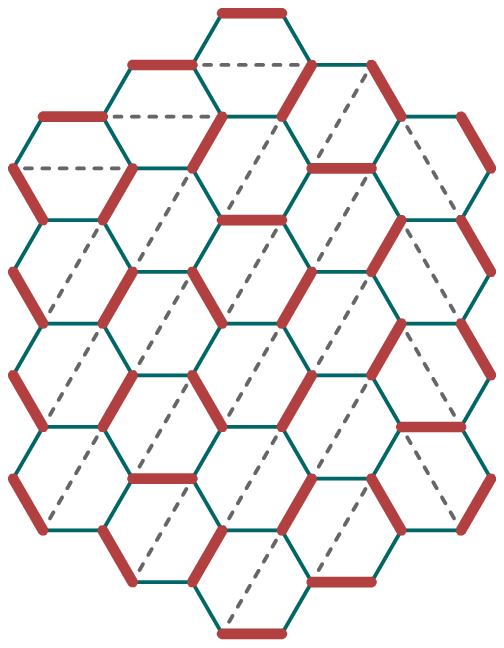}\cr
\epsfbox{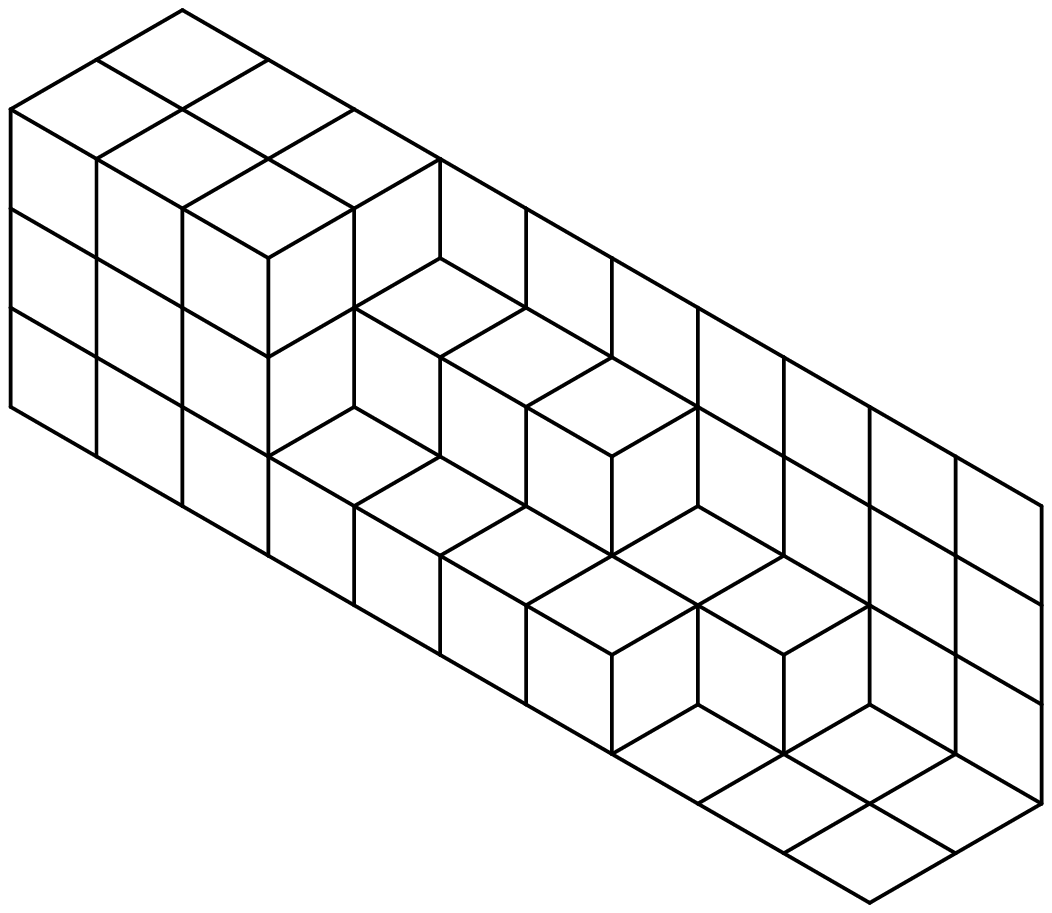}&\epsfbox{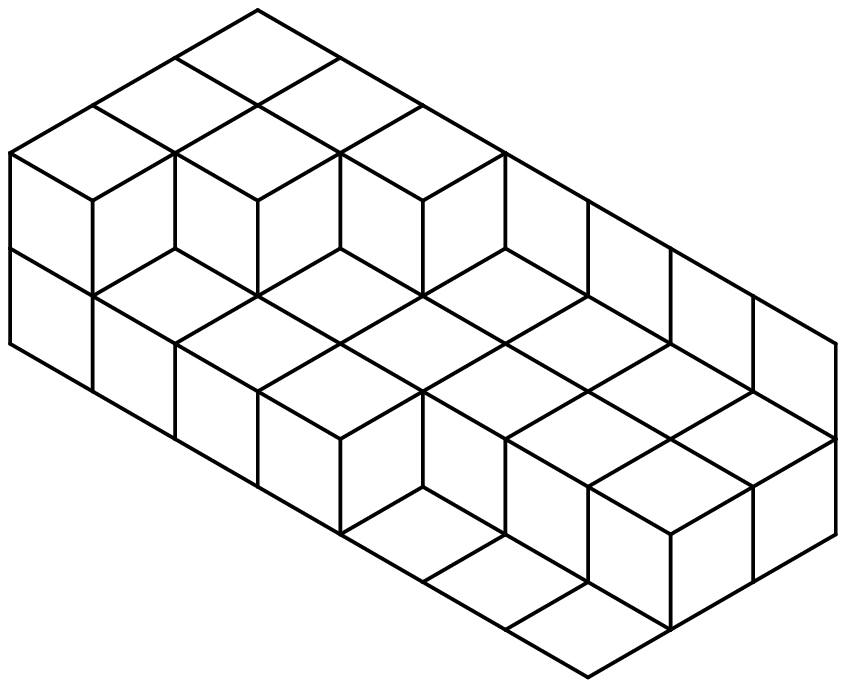}&\epsfbox{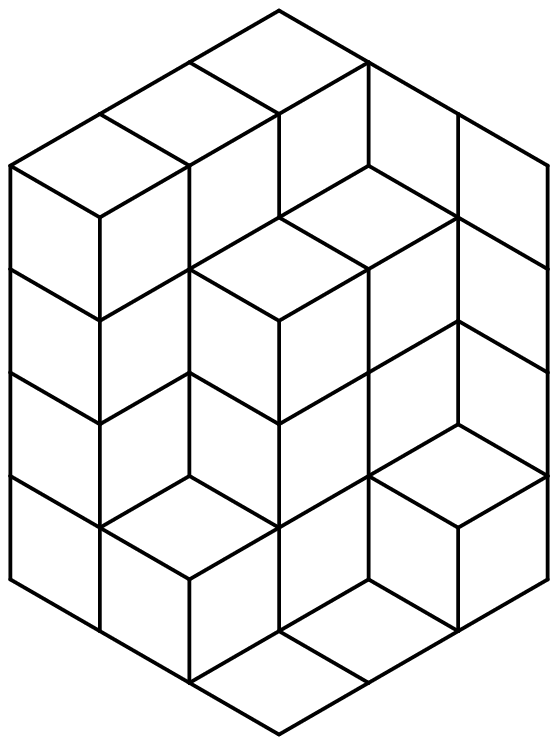}\cr
(i)&(ii)&(iii)\cr
}$
}

\newsec{Proof of the Theorem}

An important part of the proof is devoted to the determination 
of the {\it fixed edges}, i.e.\ the edges which are occupied 
(or unoccupied) in
any configuration of the given type $(A,B,C)$. We first recall
a  very useful lemma proved by de Gier \dG

\noindent{\sc Lemma 1.} {\sl In Fig.~\degier, if (i) the edges $ab$ and $ef$ 
are occupied, with $ab$ and $ef$ belonging to different loops, and 
if (ii) $cd$ either 
is an unoccupied external edge, or
belongs to a third loop, or is connected to 
$ab$ by $da$ or to $ef$ by $de$, then the edge $kl$ is occupied.}

{%
\fig\degier{de Gier's lemma}{\epsfbox{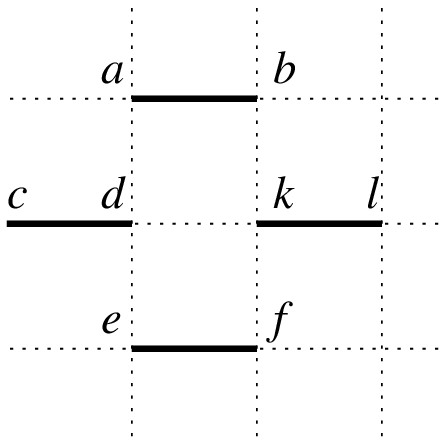}}
}

\noindent Proof of lemma 1: 
if $kl$ was not occupied, either
$bk$ and $kf$ would, which would contradict {\sl (i)}, or $dk$ and either
$kb$ or $kf$ would, which would contradict {\sl (ii)}.

\noindent{\sc Lemma 2.} {\sl The edges outside the polygon 
$P'$, as well as some edges inside the polygon, are fixed 
as depicted on Fig.~\fixedc.}

\noindent Proof of lemma 2: 
we shall build the fixed edges in 2 steps (cf Fig.~\fixeda--\fixedc).

\fig\fixeda{Empty grids with the polygon drawn.}{
$\matrix{
\epsfbox{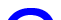}&\epsfbox{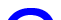}&\epsfbox{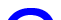}\cr
(i)&(ii)&(iii)\cr
}$ } 

\noindent
 1. We first prove that in each triangular domain limited by one of the
diagonals starting from $A$, $B$ or $C$ and the  sides of the square, 
 and exterior to the polygon $P$, the edges are fixed and form
 ``stairs''.
The external edges incident to such a triangle connect to 
external edges that are on the other side of this diagonal
and the ``loops'' that start from them
 must thus cut the diagonal at distinct points. As there are as many
points on this diagonal  as there are external
edges, the only possibility is the stair pattern. 

\fig\fixedb{Fixed edges outside the polygon,
either occupied (thick lines)
or unoccupied (thin lines).}{
$\matrix{
\epsfbox{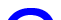}&\epsfbox{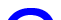}&\epsfbox{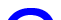}\cr
(i)&(ii)&(iii)\cr
}$
}

\noindent
2. Edges  inside the polygon may also be fixed by repeated action 
of lemma 1, starting from the sides of the square.
The fact that all the external edges 
that lie on the same side of any of the diagonals belong to different 
loops enables one to iterate the application of  lemma 1.
While this lemma fixes every other edge in a given 
direction (horizontal or vertical), {\it  all}\/ horizontal edges
in the rectangle $CB'DA''$ are fixed, by a successive application of the 
argument starting from the left and from the right external edges.

Note that at this stage, each vertex of the square grid belongs 
to at least one fixed edge.  
At the vertices belonging to two occupied edges, the 
complementary edges are also fixed to be unoccupied
 (thin solid lines on 
Fig.~\fixedc).
Those vertices  which belong to only 
one internal fixed edge may be regarded as active, 
two or three unfixed lines emanate from them 
(dashed lines on Fig.~\fixedc) and it could be
possible to switch to a dual picture, by depicting their dual triangles 
and looking at the various ways one may assemble them into lozenges, etc. 
This is the route followed by de Gier \dG\ and Krattenthaler \K. 

\fig\fixedc{The full set of fixed edges.}{
$\matrix{
\epsfbox{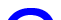}&\epsfbox{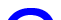}&\epsfbox{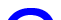}\cr
(i)&(ii)&(iii)\cr
}$
}

\noindent{\sc Lemma 3.} 
{\sl There exists a FPL configuration of type $(A,B,C)$. It is loopless.}

\noindent Proof of lemma 3: 
one may construct explicitly the two special FPL
configurations which are mapped onto the empty or the full plane
partition.
Here we shall describe the construction for one of the two, and 
we let the reader repeat it for the other. For
this purpose, we carry out a further splitting of the domains of our 
grid, and construct the points $T$,  $E$, $H$, $H'$ and $H''$ 
(see Fig.~\split). In each of the
domains limited by the polygon $P'$ and by dashed lines, 
the yet unfixed edges are determined  
according to the indicated  prescription (in red). 
It is now easy to see that  the $a$ loops entering the grid
between $A$ and $H$ will exit through $AH'$, and likewise 
the $b$ loops entering
in $BH$ exit through $BH''$ and the $c$ loops entering in $CH''$ exit
through $CH'$. This is illustrated on Fig.~\split, in which  
examples of staircase loops are represented
in each of the domains bounded by the dotted (green) lines
(see also Fig.~\empty\ for the pictures in our 
three typical examples). We thus have constructed explicitly a
configuration of type $(A,B,C)$, and it has no closed loop. 
Note that the point $E$ 
 in our construction has an interpretation 
in the FPL $\leftrightarrow$ plane partition matching: it corresponds
to the corner of the (say empty) box.

{%
\fig\split{FPL configurations corresponding to the ``empty'' plane
partition: the generic case.}{\epsfbox{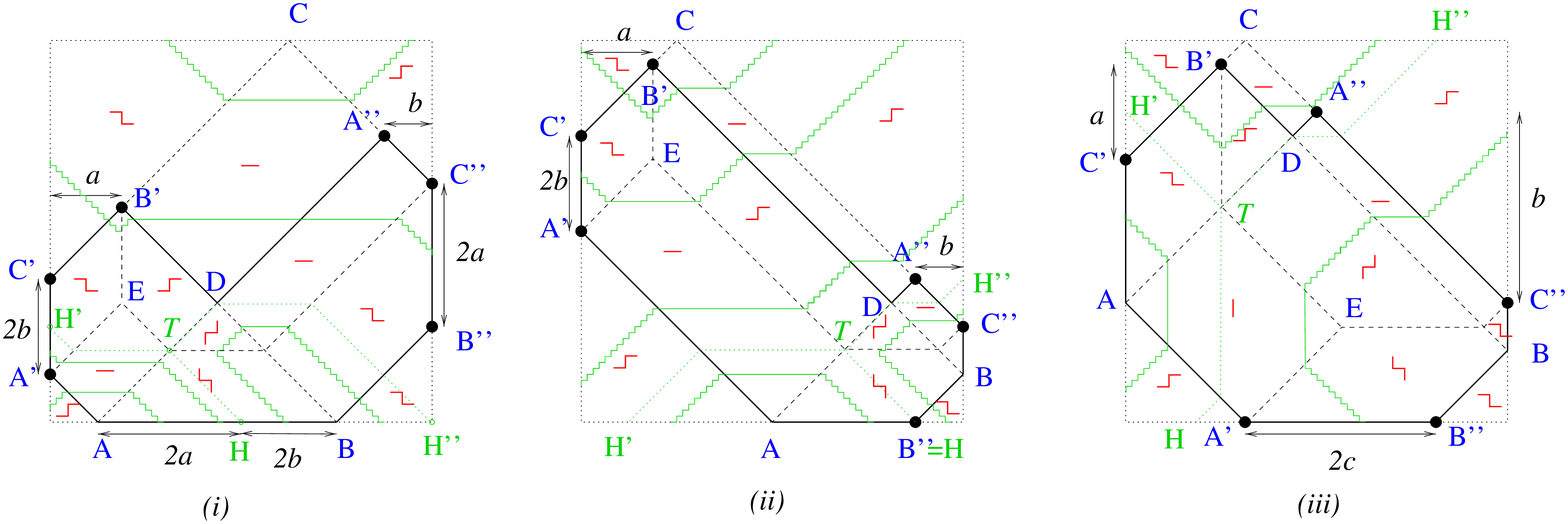}}}

\fig\empty{FPL configurations corresponding to the ``empty'' plane
partition in our three examples.}{
$\matrix{
\epsfbox{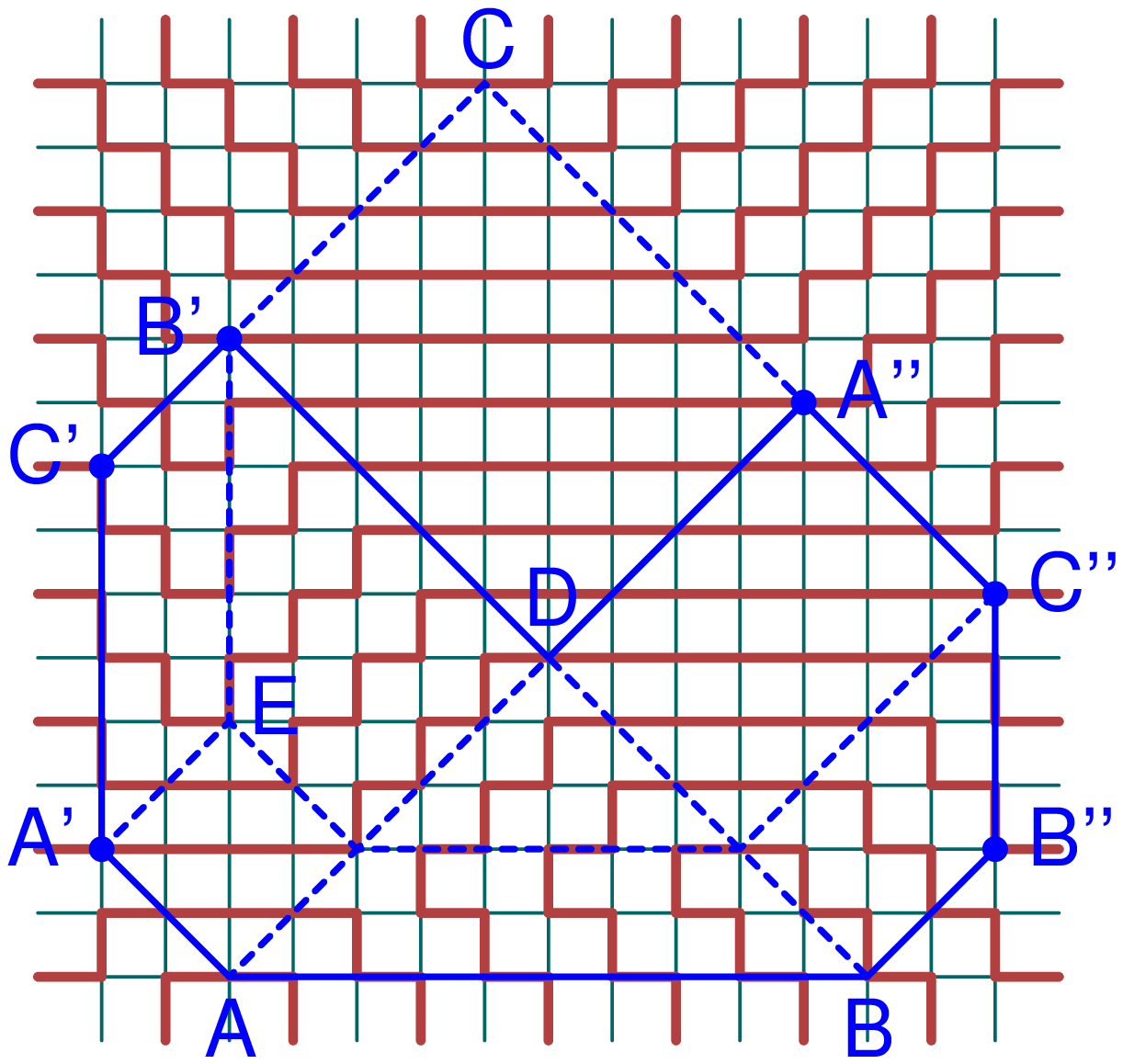}&\epsfbox{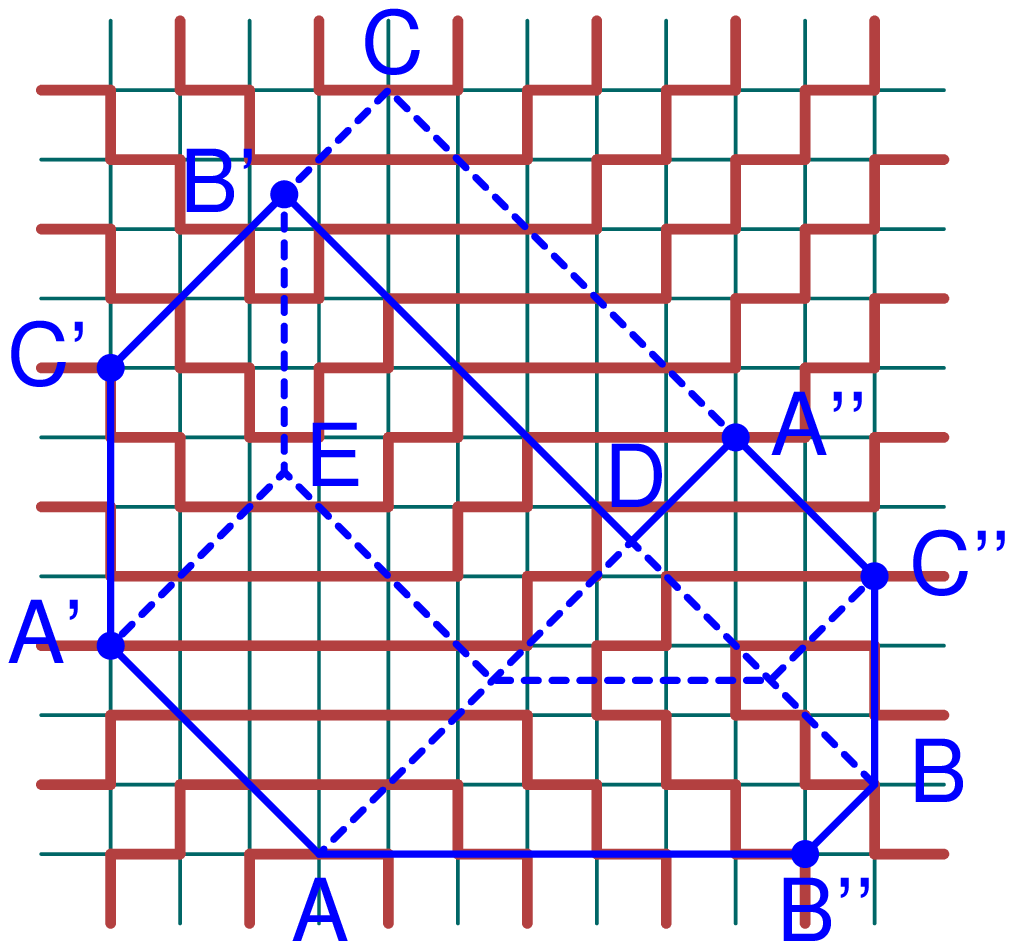}&\epsfbox{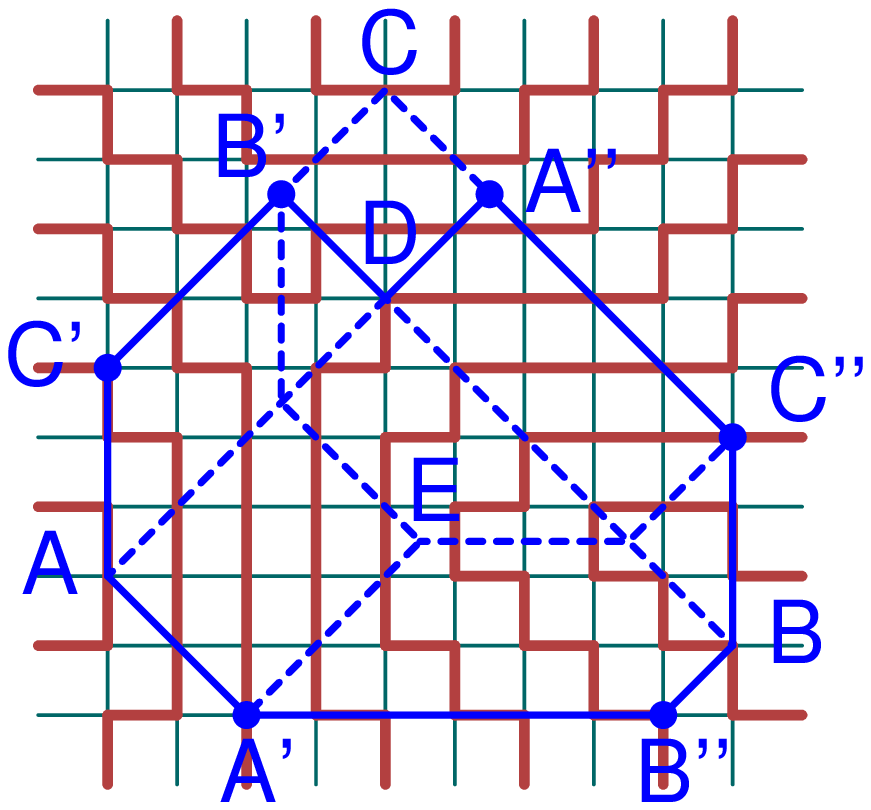}\cr
(i)&(ii)&(iii)\cr
}$
}

\medskip
\noindent Proof of theorem: 
we now have all the elements to complete the proof of the theorem.
Let us go back to the construction of section 2, step 4. To a FPL of
type $(A,B,C)$ we associate the configuration obtained by
restriction to the dominos inside the polygon $P'$. Note that
the edges of the dominos are exactly the complementary 
set of the fixed edges of Fig.~\fixedc; therefore lemma 2 implies that
this mapping is injective. By simple inspection one can check
that in all three cases depicted on Fig.~\polygon, the dominos can
be deformed into an $(a,b,c)$ hexagon. Furthermore, since exactly
two edges around each vertex are occupied on the original square
lattice,  and the
middle edges of the dominos are fixed to be occupied, the new configuration
on the hexagonal lattices is precisely a dimer configuration, or in the
dual language, a tiling of the $(a,b,c)$ hexagon with lozenges,
or still equivalently, a plane partition in a box of size $a\times b\times c$.

We thus have an injective mapping from FPL of type $(A,B,C)$ into plane
partitions. To prove its surjectivity, we note the following. The
moves  \ppmove\ on elementary hexagons of plane
partitions, or equivalently the moves \dimermove\ on dimer
configurations of the
hexagonal lattice are well-known to be ergodic,
namely allow to explore the full set of plane 
partitions in an $a\times b\times c$ box. In the 
correspondence
above, an elementary hexagon becomes deformed into a domino, and once
the fixed middle edge is added, 
the move becomes \fplmoveh\ or \fplmovev. 
A key remark is that this move does not modify the connection between 
the 4 corners of the domino; hence, it preserves
the link pattern of the whole configuration (and does not add or
remove any loops).
Starting from a particular configuration, which we choose to be the
one  exhibited in lemma 3, one
can produce using such moves the preimage of any plane partition. 
Therefore the mapping is surjective.


Alongside this proof, we have established two corollaries:

\noindent{\sc Corollary 1.} {\sl The moves \fplmovev\ and
 \fplmoveh\ are 
ergodic on FPL configurations of type $(A,B,C)$.}

\noindent{\sc Corollary 2.} {\sl There are no internal loops in any
FPL  configuration of type $(A,B,C)$.}

Both properties had been stated by de Gier \dG\ but no 
detailed proof had been given.


\newsec{Concluding comments}

The reader who is not yet fully convinced is invited to 
visit the site\hfill\break
\centerline{\tt
http://ipnweb.in2p3.fr/lptms/membres/pzinn/fpl}
to practice a little and to enjoy the show! 

Since plane partitions come with a natural grading (the number of
boxes) (up to the ambiguity with the complement), 
our bijection  yields a grading of FPL of type $(A,B,C)$
and Macdonald formula
$ \prod_{i=1}^a\prod_{j=1}^b\prod_{k=1}^c {1-q^{i+j+k-1}\over 1-q^{ i+j+k-2}} $
 gives the relevant counting. It would be 
nice to have a direct interpretation of this grading in the language of FPL, 
and to see if it extends to other types of FPLs. 

By Wieland's theorem \Wie, it is known that the number of FPL of
a given link pattern is invariant under the action of the dihedral 
group $D_{2n}$ on this link pattern. 
One observes empirically that the effect of a Wieland's rotation 
is to rotate the diameters of the elementary hexagons of step~3 in
sec.~2, see also Fig.~\polygon,  
while the dimer configuration is left 
invariant. In other words, a FPL configuration and its Wieland rotation
share the same image under the respective bijections. 
Note that Wieland's rotations do not create any loop in this case, 
even though they do in general. 

Note that the two problems connected by our bijection, 
namely  the FPL counting and  the plane 
partition counting, may both be rephrased in terms of the 6-vertex model.
They correspond however to different weights and to different 
fixed boundary conditions. 

\fig\fulempty{Full/empty configurations in the HFPL picture,
for $(a,b,c)=(5,4,3)$
.}{\epsfbox{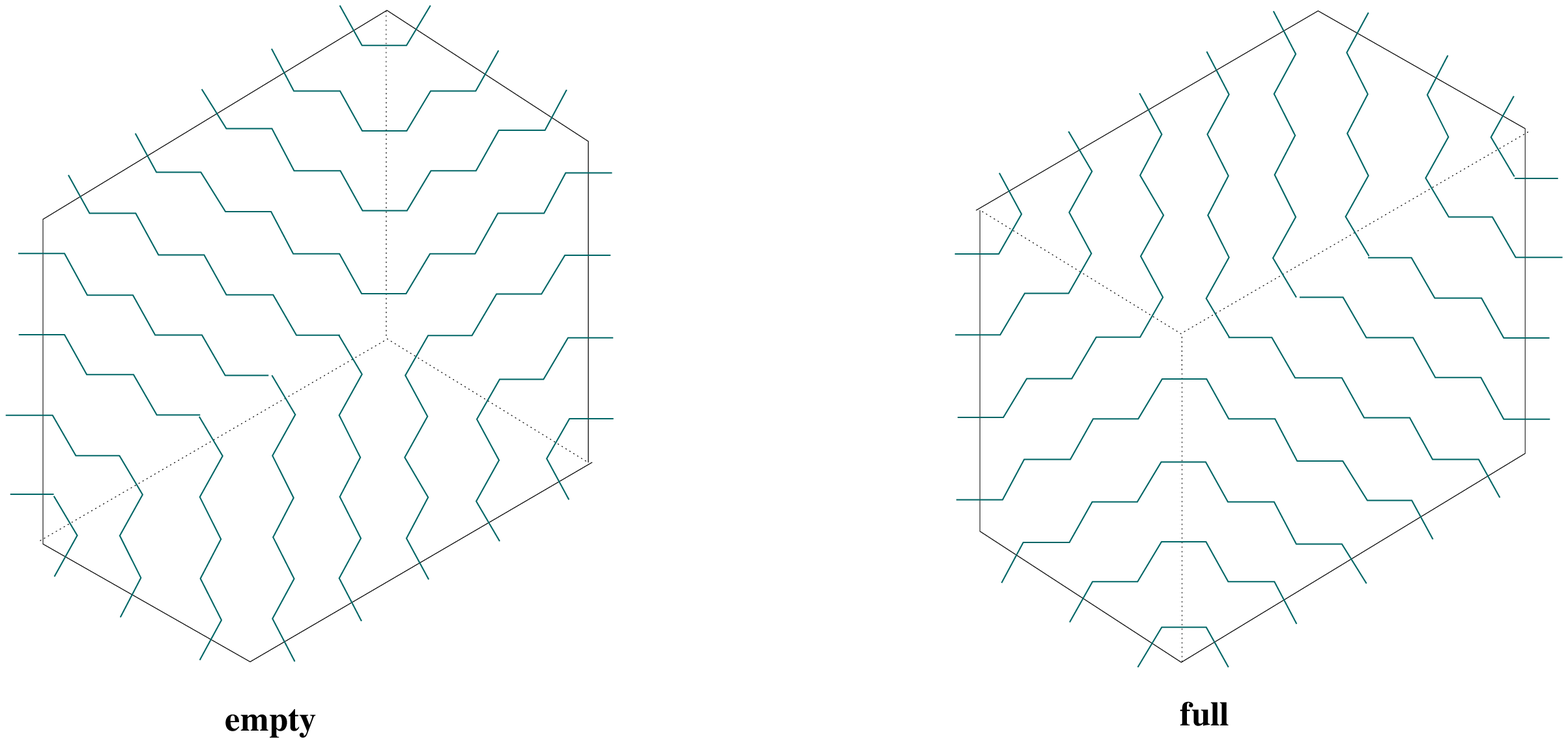}}

A last remark is in order. The dimer formulation of
rhombus tilings or plane partitions in an $a\times b\times c$ box
may also be rephrased by considering the ``loops'' formed by the
edges complementary to the dimers on the underlying
hexagonal lattice.  These loops are paths that 
connect points on the perimeter of the $(a,b,c)$ hexagon
and form Fully Packed Loop
configurations of the hexagonal lattice, referred to as HFPL. 
The HFPL are nothing but the images of the complement
of the FPL of type $(A,B,C)$ in our bijection.
As is readily seen, 
these complementary loops connect every other point on
the perimeter of the square to every point on the perimeter
of the domino tiling of the active zone (i.e.\ the polygon $P'$), 
via ``parallel'' paths.
Thus, our bijection also yields a one-to-one mapping of
the FPL's complement to the HFPL, and in particular the link
patterns of the latter yield those of the former. Among the HFPL
two are particularly simple: they correspond to the two fundamental
states of an empty or a completely filled plane partition, see
Fig.~\fulempty. As the HFPL paths must travel along parallel 
zig-zag lines on the three
visible sides of the (empty or filled) box,  their
link pattern is also of the type $(a,b,c)$, and so is that of the
two corresponding FPL's complements, which therefore have no
internal loop. The complement of a generic FPL of type $(a,b,c)$, however,
may have internal loops, as the basic move does not preserve the connectivity
of the complementary edges. The identification with HFPL
provides therefore yet another grading of FPL, according to the number of
internal loops of the complement or equivalently of the associated HFPL.
Conversely, this allows to view all HFPL of an $(a,b,c)$ hexagon
as the complement of the restriction of FPL to $(a,b,c)$ link patterns.
One may wonder whether more general FPL's complements could correspond to
HFPL on more general domains.

{\vskip1cm
\centerline{\bf Acknowledgments} 
It is a pleasure to thank J.~de Gier, J.~Propp and D. Wilson
for stimulating comments and especially C.~Krattenthaler for
informing us of his own work on proving the bijection and for 
fruitful exchanges. 
This work is partially supported by the European networks  
HPRN-CT-2002-00325 and HPRN-CT-1999-00161.

 }

\listrefs
\bye